%% file: main.tex
\pdfoutput=1
\documentclass[conference]{IEEEtran}
\IEEEoverridecommandlockouts
\usepackage{cite}
\usepackage{amsmath,amssymb,amsfonts}
\usepackage{algorithmic}
\usepackage{graphicx}
\usepackage{textcomp}
\usepackage{xcolor}

\usepackage{multirow}
\usepackage{url}
\usepackage[utf8]{inputenc}
\usepackage{cleveref}

\usepackage{nomencl}
\usepackage{tabularx,booktabs}
\usepackage[linesnumbered]{algorithm2e}

\def\BibTeX{{\rm B\kern-.05em{\sc i\kern-.025em b}\kern-.08em
    T\kern-.1667em\lower.7ex\hbox{E}\kern-.125emX}}

\makeatletter
\def\ps@IEEEtitlepagestyle{%
  \def\@oddfoot{\mycopyrightnotice}%
  \def\@evenfoot{}%
}
\def\mycopyrightnotice{%
  {\hfill \footnotesize 
  
  reprinting/republishing this material for advertising or promotional purposes,creating new collective works, for resale or redistribution to servers or lists, or reuse of any copyrighted component of this work in other works. \hfill}
}
\makeatother

\makeatletter
\def\ps@IEEEtitlepagestyle{
  \def\@oddfoot{\mycopyrightnotice}
  \def\@evenfoot{}
}
\def\mycopyrightnotice{
  {\footnotesize
  \begin{minipage}{\textwidth}
  ©2020 IEEE Personal use of this material is permitted. Permission from IEEE must be obtained for all other uses, in any current or future media, including
reprinting/republishing this material for advertising or promotional purposes,creating new collective works, for resale or redistribution to servers or lists, or
reuse of any copyrighted component of this work in other works.

  \end{minipage}
  }
}

\begin{document}

\title{Representing Long-term Impact of Residential Building Energy Management using Stochastic Dynamic Programming \\
%Using stochastic dynamic programming to represent long-term impact of building scheduling\\
%Long-term optimal scheduling of a residential building with a capacity-based grid tariff when considering uncertainty
\thanks{Funded by FME CINELDI/FME ZEN/Norwegian Research Council, C.no 257626/257660/E20}
%}
% \author{\IEEEauthorblockN{1\textsuperscript{st} Kasper Emil Thorvaldsen}
% \IEEEauthorblockA{\textit{Department of Electric Power Engineering} \\
% \textit{Norwegian University of Science and Technology (NTNU)}\\
% Trondheim, Norway \\
% kasper.e.thorvaldsen@ntnu.no}
% \and
% \IEEEauthorblockN{2\textsuperscript{nd} Sigurd Bjarghov, 3\textsuperscript{rd} Hossein Farahmand}
% \IEEEauthorblockA{\textit{Department of Electric Power Engineering} \\
% \textit{Norwegian University of Science and Technology (NTNU)}\\
% Trondheim, Norway \\
% email address or ORCID}
% %\and
% %\IEEEauthorblockN{3\textsuperscript{rd} Hossein Farahmand}
% %\IEEEauthorblockA{\textit{Department of Electric Power Engineering} \\
% %\textit{Norwegian University of Science and Technology (NTNU)}\\
% %Trondheim, Norway \\
% %email address or ORCID}
% 
}
\author{\IEEEauthorblockN{Kasper Emil Thorvaldsen*\thanks{*Corresponding author.}, Sigurd Bjarghov, Hossein Farahmand}
\IEEEauthorblockA{\textit{Department of Electric Power Engineering} \\
\textit{Norwegian University of Science and Technology (NTNU)}\\
Trondheim, Norway \\
kasper.e.thorvaldsen@ntnu.no, sigurd.bjarghov@ntnu.no, hossein.farahmand@ntnu.no }}

% Tips fra meg: Da får man samlet alle til en

%\author{\IEEEauthorblockN{Kasper Emil Thorvaldsen, Sigurd Bjarghos, and Hossein Farahmand}
%\IEEEauthorblockA{\textit{Department of Electric Power Engineering} \\
%\textit{Norwegian University of Science and Technology (NTNU)}\\
%Trondheim, Norway \\
%\{kasper.e.thorvaldsen@ntnu.no, email2, and email3 \}}}

\maketitle

\begin{abstract}
Scheduling a residential building short-term to optimize the electricity bill can be difficult with the inclusion of capacity-based grid tariffs. Scheduling the building based on a proposed measured-peak (MP) grid tariff, which is a cost based on the highest peak power over a period, requires the user to consider the impact the current decision-making has in the future. Therefore, the authors propose a mathematical model using stochastic dynamic programming (SDP) that tries to represent the long-term impact of current decision-making. The SDP algorithm calculates non-linear expected future cost curves (EFCC) for the building based on the peak power backwards for each day over a month. The uncertainty in load demand and weather are considered using a discrete Markov chain setup. %, and the operational scheduling for a day to find the optimal peak is influenced by the computed EFCC.
The model is applied to a case study for a Norwegian building with smart control of flexible loads, and compared against methods where the MP grid tariff is not accurately represented, and where the user has perfect information of the whole month. The results showed that the SDP algorithm performs 0.3 \% better than a scenario with no accurate way of presenting future impacts, and performs 3.6 \% worse compared to a scenario where the user had perfect information. %, and performs 4.2 \% worse than if the user had perfect information. The results also captured areas of improvement that can give more accurate representation of the future.
%Still, the algorithm captures the benefits of representing the future and located areas of improvement for further work on this concept.

\end{abstract}


\input{chapters/Nomenclature}
\printnomenclature

\input{chapters/Introduction}

%\section{Methodology}

\input{chapters/Methodology}
\input{chapters/Case_study}
\input{chapters/Results_and_discussions}
\input{chapters/Conclusion}

\bibliographystyle{IEEEtran}
\bibliography{IEEEabrv,references}

\end{document}

%% file: chapters/Nomenclature.tex
\makenomenclature

\renewcommand{\nomgroup}[1]{%
\ifthenelse{\equal{#1}{A}}{\item[\textbf{SDP sets}]}{%
\ifthenelse{\equal{#1}{B}}{\item[\textbf{Index sets}]}{%
\ifthenelse{\equal{#1}{C}}{\item[\textbf{Parameters}]}{%
\ifthenelse{\equal{#1}{D}}{\item[\textbf{Decision variables}]}{%
\ifthenelse{\equal{#1}{E}}{\item[\textbf{Stochastic variables}]}{%
\ifthenelse{\equal{#1}{O}}{\item[\textbf{Other}]}{}}}}}}
}

%missing nomenclature:
%Question: Should SOS-2 be in the nomenclature, or is it enough with the explanation added? Does the same go for "Stochvar"?

%C^p_n
%stochvar
%probability

%For sets
%\nomenclature[A]{$\mathcal{T}_{Stage}$}{set of stages within the period}
\nomenclature[B]{$\mathcal{S}_{g}$}{set of state variables}
%\nomenclature[BA]{$\mathcal{S}_{s,g}$}{set of stochastic variables}
\nomenclature[BA]{$G$}{set of days within the month}
%\nomenclature[B]{$\mathcal{S}_{p,\tau}$}{set of deciding variables}
%\nomenclature[A]{$\mathcal{N}_{State}$}{List of discrete}
%\nomenclature[A]{$\mathcal{N}_{Scen}$}{set of deciding variables}

\nomenclature[B]{$\mathcal{T}$}{Set of time steps within the day}

%Stochastic variables
\nomenclature[E]{$P^{spot}_t$}{Electricity spot price in time step t [$\frac{EUR}{kWh}$]}
\nomenclature[E]{$D_{t}^{el}$}{Consumer-specific load in time step t [$kWh$]}
\nomenclature[E]{$\delta_t^{EV}$}{EV connected to building $\{0,1\}$}
\nomenclature[E]{$y^{PV}_t$}{Photovoltaic production from installed system [$kWh$]}
\nomenclature[E]{$\delta_t^{Occ}$}{Occupancy presence $\{0,1\}$}
\nomenclature[E]{$T^{out}_t$}{Daily average outdoor temperature [$^\circ C$]}
%\nomenclature[E]{$D^{DHW}_t$}{Domestic hot water load [$kWh$]}

%Decision variables
\nomenclature[D]{$y_t^{imp}, y_t^{exp}$}{Energy imported/exported to household [$\frac{kWh}{h}$]}
%\nomenclature[D]{$p_g^{exp}$}{Energy exported from household [$\frac{kWh}{h}$]}
\nomenclature[D]{$p^{p}$}{Peak of imported energy [$\frac{kWh}{h}$]}
\nomenclature[D]{$\gamma$}{SOS-2 variables for the Expected future cost curve}
%Special variable
\nomenclature[D]{$C_{y^{p},s^e_{t+1}}^{p}$}{Expected future cost from peak power [$\frac{EUR}{\frac{kWh}{h}}$]}

%EV
\nomenclature[D]{$SoC_t^{EV}$}{State of charge for EV for time step t [$kWh$]}
\nomenclature[D]{$y_t^{EV,ch}$}{Input power to EV for time step t [$\frac{kWh}{h}$]}

%BAT
\nomenclature[D]{$SoC_t^{B}$}{State of charge for Battery for time step t [$kWh$]}
\nomenclature[D]{$y_t^{B,ch},y_t^{B,dch}$}{Power to/from the battery for time step t [$\frac{kWh}{h}$]}

%PV
%\nomenclature[D]{$y_t^{PV}$}{Power produced from PV system [$\frac{kWh}{h}$]}

%Heat
\nomenclature[D]{$T_t^{in},T_t^{e}$}{Interior and building envelope temperature [$^\circ C$]}
%\nomenclature[D]{$T_t^{e}$}{Building envelope temperature [$\degree C$]}
%\nomenclature[D]{$T_t^{w}$}{Water tank temperature [$\degree C$]}

\nomenclature[D]{$q_t^{sh}$}{Power usage for space heating [$\frac{kWh}{h}$]}
%\nomenclature[D]{$\mu_g^{\Delta T}$}{Change in energy in room [$\frac{kWh}{h}$]}

%Parameters

\nomenclature[C]{$C^{grid}$}{DSO energy tariff for imported energy [$\frac{EUR}{kWh}$]}
\nomenclature[C]{$P^{imp,max}$}{Maximum power import to building [$\frac{kWh}{h}$]}
\nomenclature[C]{$P_n^{p}$}{Peak power at point $n$ [$\frac{kWh}{h}$]}
\nomenclature[C]{$P_0^{p}$}{Initial peak power [$\frac{kWh}{h}$]}
\nomenclature[C]{$VAT$}{Value added tax for purchase of electricity [$p.u$]}

%EV
\nomenclature[C]{$E^{EV,Cap}$}{EV storage capacity [$kWh$]}
\nomenclature[C]{$\eta^{EV}_{ch}$}{EV charging efficiency [$\%$]}
\nomenclature[C]{$\dot{E}^{Max}$}{Maximum EV charging capacity [$\frac{kWh}{h}$]}
\nomenclature[C]{$SoC^{EV,min},SoC^{EV,max}$}{Min/Max EV soc capacity [$kWh$]}
\nomenclature[C]{$D^{EV}$}{Hourly EV discharge when not connected [$kWh$]}
%\nomenclature[C]{$SoC^{EV,Start/Stop}$}{EV start/end soc [$kWh$]}

%HERE I AM, ONCE AGAIN!
%BAT
\nomenclature[C]{$E^{B,Cap}$}{Battery storage capacity [$kWh$]}
\nomenclature[C]{$\eta^{B}_{dch},\eta^{B}_{ch}$}{Discharge/charge efficiency for battery [$\%$]}
\nomenclature[C]{$\dot{E}^{B,dch},\dot{E}^{B,ch}$}{Discharge/charge capacity for battery [$\frac{kWh}{h}$]}
%\nomenclature[C]{$SoC^{B,Start/Stop}$}{Battery start/end soc [$kWh$]}
\nomenclature[C]{$SoC^{B,min},SoC^{B,max}$}{Battery soc limits [$kWh$]}

%PV
%\nomenclature[C]{$A^{PV}$}{PV system area [$m^2$]}
%\nomenclature[C]{$\eta^{PV}$}{Total efficiency for PV system [$\%$]}

%Heat
\nomenclature[C]{$\dot{Q}^{sh}$}{Capacity for space heating radiator [$\frac{kWh}{h}$]}
%\nomenclature[C]{$T^{in,Start/End}$}{Start/end temperature of the building [$\degree C$]}
\nomenclature[C]{$T_t^{in,min}, T_t^{in,max}$}{Lower/upper interior boundary [$^\circ C$]}
%\nomenclature[C]{$C^{building}$}{specific heat capacity for heating storage in building [$\frac{kWh}{m^2\cdot K}$]}
%\nomenclature[C]{$\alpha_g^{HL}$}{Time-dependent heat loss constant [$kWh$]}
%\nomenclature[C]{$\beta_g^{HL}$}{Time-dependent heat loss slope [$\frac{kWh}{\degree C}$]}
%\nomenclature[C]{$M^{Occ}$}{Heat gained due to occupancy [$\frac{kWh}{h}$]}
\nomenclature[C]{$R_{ie}, R_{eo}$}{The thermal resistance between the interior-building envelope and building envelope-outdoor air [$\frac{^\circ C}{kWh}$]}
\nomenclature[C]{$C_i, C_e$}{Heat capacity for interior and building envelope [$\frac{kWh}{^\circ C}$]}

%Water
%\nomenclature[C]{$\dot{Q}^{w,rated}$}{Rated capacity for water tank heater [$\frac{kWh}{h}$]}
%\nomenclature[C]{$T^{w,min}, T^{w,max}$}{Lower/upper temperature boundary for water tank[$\degree C$]}

%Nodes
\nomenclature[C]{$N_S$}{Number of nodes for stochastic variables}
\nomenclature[C]{$N_P$}{Number of discrete peak power values}
%\nomenclature[C]{$T^{amb}$}{Ambient temperature for water tank [$^\circ C$]}

%For parameters

%% file: chapters/Introduction.tex
\section{Introduction} \label{Introduction}

In Norway, the Norwegian regulator (NVE) has proposed that the grid tariff must be updated to better reflect the actual costs associated with the operating of the distribution grid \cite{NVEEnergidirektorat}. The nationwide rollout of smart meters gives an opportunity to implement more sophisticated pricing structures. Residential consumption is trending towards a more peak intensive profile, requiring investments to cover the relatively few peak-load and scarcity hours in the distributional grid. Therefore, NVE has proposed a new grid tariff with cost elements related to both capacity and energy.

%resulting in more scarcity hours in the distribution grid, requiring investments for relatively few peak-load hours. 

In a hearing from NVE, they discuss two types of capacity-based tariffs; measured-peak (MP) and subscribed capacity tariffs \cite{Norgesvassdrags-ogenergidirektorat2017ForslagNettvirksomhet}. Measured-peak prices the highest measured peak as the base for the cost the customer has to pay over a specified period, whereas subscribed capacity requires the customer to subscribe to a certain amount of capacity, and then pays a low or high energy cost if the consumption is under or over the subscribed capacity, respectively. The MP grid tariff is affected by the active import of electricity, and can be costly if the end-user is not aware of their electricity consumption. The proposed period for the MP grid tariff for residential buildings is 1 day, but this tariff already exists for many commercial/industrial customers today with periods over 1 month. Thus, it is not unlikely that this tariff could be implemented for longer periods for residential buildings, making it important to consider previous and future actions to prevent a costly electrical bill.
%If this tariff would be extended to cover a longer period for residential buildings, then one will need to be aware of the consumption and previous action, but also have an indication of the future impact of the scheduling of the building in terms of possible load outcome and weather forecast. This uncertainty in future situations makes it key to represent it for the end-user as a prediction.

With the smart meters rollout in over Norway, it is only a matter of time before smart controllers for the flexible appliances in the building are common. Smart building controllers try  to reduce the electricity bill using the optimal scheduling of available flexibility in the buildings on a short time-scale. Reference \cite{Zhao2015MPC-basedStorages} optimizes the building flexibility using a model predictive control (MPC)-based optimal scheduling strategy with a non-linear programming model for the next 24 hours. Moreover, it minimizes cost and CO2 emissions simultaneously. In \cite{BehzadiForough2017MultiSystem}, authors present a methodology for an energy management system, aiming to find optimal scheduling of a hybrid renewable energy system. The methodology is based on multi-objective receding horizon optimization, and investigated optimal scheduling using a time view from 6-24 hours with predicted inputs. The findings showcase that a longer time horizon in input predictions, contribute to increase the renewable energy share for covering the demand. Dynamic programming was used in \cite{Jones2018SolvingCharges} to operate a battery optimally. The objective is to minimize the electricity cost subject to grid tariffs based on peak consumption. The battery is optimized over 24 hours with hourly steps, and tested for both a deterministic and stochastic setup to reduce the demand charge costs.
%There were no future costs considered beyond the analyzed period. 

The above literature consider the short-term optimization, with uncertainty reflected in the hourly future impact. However, the analysis period and the number of hours/days into the future limits the amount of information they have. Mostly, there is no need to go further as short-term optimization utilizes flexibility that is shifted a couple of hours, and has no indicated impact in the future. However, if longer periods must be considered, like what grid tariffs promotes, then it is required to connect the future unknown scheduling of flexibility realization together with daily operational decision-making to make more accurate decisions.

This problem has been addressed previously regarding hydropower scheduling, where the water stored in the reservoirs should be used optimally. As presented in \cite{Mo2002Short-termMarket}, the planning problem for a hydropower producer is complex due to the interactions with power systems with thermal power production. To find optimal operation, the problem is decomposed into several smaller optimization problems, for instance long-, medium-, and short-term scheduling. The long-term scheduling includes uncertainty and the objective is to find the marginal value of using water now to produce power versus the value of storing it for later. The calculation is done using stochastic dynamic programming (SDP) to find these marginal values called water values, as also shown in \cite{Helseth2017AssessingMarkets, Stedinger1984StochasticOptimization}. The ending water values are used in a short-term model to optimize production. This setup is able to represent the future impact of current decision-making, and could therefore be implemented into a building model as well to complement the problems discussed. Using the proposed MP grid tariff is a good way of analyzing the performance of such a setup. Dynamic Programming was used in \cite{Jones2018SolvingCharges} to optimally operate a battery to minimize cost under MP grid tariffs, but only considers a 24 hour horizon, which would be insufficient with a monthly tariff where uncertainty for 1 month must be taken into consideration.

In this paper, we propose a backwards SDP algorithm that generates non-linear future cost curves for a residential building for every day over the course of a month, inspired by the long-term model used frequently in hydropower scheduling. Our contributions are the following:
\begin{itemize}
    \item We present a SDP optimization framework to minimize costs by predicting and planning for the highest peak with uncertain demand
    \item We compare the results of the SDP algorithm to cases where the MP grid tariff is either considered short-term or not considered, and with cases where the residential building has perfect information for the entire month
\end{itemize}

The rest of the paper is organized as follows: \Cref{Model} describes the mathematical formulation and the SDP layout. \Cref{Case_study} describes the performed case study, whereas \Cref{discussion} presents and discusses the results. The conclusion is found in \Cref{Conclusion}.

%% file: chapters/Methodology.tex
\section{Methodology} \label{Model}

The methodology takes into consideration a single all-electric residential building that is connected to the power grid with bi-directional power flow options. The building considers smart control of different application such as space heating, electric vehicle (EV) charging, battery control and a photovoltaic (PV) system connected on the roof. %The smart control is owned by the end-user, which will be the participant benefiting from this methodology.

\subsection{Model overview}
The overall objective is to find the operation strategy for the end-user over a month that minimizes the total electricity bill for the building when a one-time MP grid tariff cost based on the highest peak power imported is included. The whole month is modelled as a multi-stage stochastic optimization problem, where the stochastic variables within the building are realized for each scheduling day with scenarios. The probability distributions for each stochastic scenario are assumed to be discrete, and it is assumed that we can decompose the problem into daily decision stages. For each day, the only information carried over are the state variables, which in this work is the achieved peak power $p^p$.

%The problem is decomposed into daily decision stages using stochastic dynamic programming (SDP), where the realization of the stochastic variables for the given day is known at the start of each day. 

The motivation behind decomposing the problem is due to the complexity of the original model, which can be affected by dimensionality issues if all possible combinations of outcomes were included. %if all the possible combinations would be included as a result of adding uncertainty, 
%For instance, with the number of scenarios per day being $n$, the total number of possible combinations over the month are $n^{\#_{days}}$. Thus, decomposition techniques to shorten down the computation time, while keeping as much information about the scenarios and system as possible is key to make this manageable.
The SDP keeps sufficient levels of information present, while still being capable of showcasing how system costs are impacted by state variables under uncertainty.
%SDP is an approach that should be capable of keeping the level of information present while decomposing to make it easier to solve, to avoid running all combined scenarios. Also, by using SDP, one is capable of showcasing how the state variables in the system is impacting the future cost for the system, even when including uncertainty, as a weighted outcome.

Some of the shortcomings of SDP, however, is that some information is lost when decomposing the problem. Only the state variables can be carried over between days, which means other variables with information must be simplified to a fixed start/end parameter value at the start/end of each day.
%Between the stages, only the information that is investigated with SDP can be transferred between, meaning that other behaviour must be simplified to be a fixed parameter value at the start and end of the decision stages.
Fixing variables during the transition will lead to loss of accuracy as the interaction between could provide better performance. The other drawback is scalability, in that if the case is complex and the number of decision stages that must be run are high, the computation time might take too long. One solution to this is to decrease the number of state variables for the SDP to calculate, giving an accuracy versus computation time issue.

When decomposing the SDP problem, we include a set of state variables $\mathcal{S}_g$ that contains all information carried over between the decision stages $g-1$ to $g$. This set is divided into two subsets. Subset $\mathcal{S}_{P,g} \in \mathcal{S}_g$ consist of the state variables that are tied together with the optimization problem, which will be explained further in Section \ref{Model_build}. Subset $\mathcal{S}_{S,g} \in \mathcal{S}_g$ contains all the stochastic variables that are being realized as parameters at the beginning of each day $g$, described in more detail in Section \ref{Stochastic_var}. Then, the decomposed decision problem for a given state $ {s_g^s,s_g^p} \in \mathcal{S}_g$ at the start of decision stage $g$ will be based on the current scheduling and the weighted impact of the future cost for all scenarios, which is the objective function of the optimization problem found in (\ref{eq:Obj_function}).

%Les mer opp på dette temaet, kanskje det kan forenkles mer!
%\textbf{The SDP formulation decomposes the scheduling period into decision stages. Then, we define a set of system states $\mathcal{S}_{t}$ that contains all information that is transferred from one decision stage $t\textendash1$ to the next stage $t$. This set is divided into two subsets; the deciding variables in $\mathcal{S}_{p,t}$ capture the consequence of the future based on the scheduling in the current stage, which is described in Section \ref{Model_build}. The stochastic variables in $\mathcal{S}_{s,t}$ contain the variables being realized in stage $t$, scenario $s$, explained further in Section \ref{Stochastic_var}. Then, the decomposed decision problem for a given state $ {s_t^e,s_t^p} \in \mathcal{S}_t$ at the start of decision stage $t$ will be based on the current scheduling and the weighted impact of the future cost for all scenarios given the state $n$. This is shown in Figure \ref{fig:EFC_behvaior}.}
%\textbf{Maybe this should be directed to Algorithm 1. This figure is as people said not very well connected. Check reviews to look for suggestions!}
%\begin{figure}[h!]
%\centering
%\includegraphics[width=0.26\textwidth]{Figs/SDP_behavior.PNG}
%\caption{Illustration of how the future cost for stage $t+1$ is represented in stage $t$ for a given state $n$ with $m$ scenarios.}
%\label{fig:EFC_behvaior}
%\end{figure}

\subsection{Stochastic variables} \label{Stochastic_var}

The stochastic variables, which are unknown values for the system until a scenario has realized their values, are assumed to be Markovian. The future cost in the optimization problem in (\ref{eq:Obj_function}) can then be represented as a Markov decision process. The Markovian setup defines the future probable states that can occur to only be dependent on the current state. Hence, we can connect the impact of the stochastic variables between the stages, and represent them for stage $g$ as an expected future cost (EFC) for all future states in stage $g+1$. The EFC will be the weighted possible outcome of all discrete scenarios $N_{S}$ in $\mathcal{S}_{S,g+1}$. In this work, there are 6 different stochastic variables per scenario: Electric-specific non-shiftable demand, electricity spot price,  EV availability, PV production, occupancy presence in the building and outdoor temperature.

\subsection{Decision problem and SDP algorithm} \label{Model_build}
The decomposed decision problem for one day $g$ is formulated as a MILP problem given in (\ref{eq:Obj_function})\textendash (\ref{eq:Envelope_balance}). The stochastic variables in $s_g^s$ is known and has been realized, and the state variable $s_g^p = (P_0^{p})$ as well.

%Note, either write it as forall t, or for t in T. Find out on monday from Christian
%For all: \quad \forall t
%For t in T: \quad t \in \mathcal{T}

%\setlength\parskip{-2em plus 0.1em minus 0.2em}

\subsubsection{Objective function}
The objective function in (\ref{eq:Obj_function}) tries to minimize the total cost of the end-user. This is denoted by the cost of purchasing electricity from the grid ($C^{Import}$), the benefit of selling excess electricity to the grid ($C^{Export}$), and the EFC based on the achieved peak power $p^p$ ($C^{Future}$).
\begin{align}
    & min \; [C^{Import} - C^{Export} + C^{Future}] \label{eq:Obj_function}\\
    & C^{Import} = \sum_{t \in \mathcal{T} } y_t^{imp}\cdot (C^{grid} + P_t^{spot})\cdot (1+VAT) \nonumber\\
    & C^{Export} = \sum_{t \in \mathcal{T} } y_t^{exp}\cdot P_t^{Spot} \nonumber \\
    & C^{Future} =  C_{p^{p},s_{t+1}^{s}}^{p}   \nonumber 
\end{align}
\subsubsection{Expected future cost}
The EFC is depicted within (\ref{eq:Memory_power}) to (\ref{eq:SOS2_var}). The highest amount of power that is imported to the building is found in $p^{p}$, which depends on the highest peak within the decision stage and the initial value given from the state variable. The peak power achieved at the end is used to set the EFC, which consist of discretized peak power of $n = 1 ... N_{P}$ points, from [$P_1^{p} = 0, P_{N_{P}}^{p} = P^{imp,max}$]. The EFC curve (EFCC) is then given as a piecewise-linear function using a SOS-2 set. SOS-2 is an ordered set of non-negative variables, where at most two variables can be non-zero, under the requirement that they are adjacent to each other in the set. The variables must sum up to 1. The SOS-2 set depicts where in the piecewise-linear function the peak power is at (\ref{eq:SOS2_binding}), and finds the resulting cost in (\ref{eq:SOS2_decision}) \cite{Beale1976GLOBALSETS}. The obtained future cost from $p^{p}$ is included into the objective function.
%e optimization problem will include this cost in the objective function.
Calculation of the EFCC is explained further in Section \ref{sol_strat}.
\begin{subequations}
    \begin{flalign}
        & p^{p} \geq P_0^{p} \label{eq:Memory_power}\\
        & p^{p} \geq y_t^{imp} \quad \forall t \label{eq:hourly_power}\\
        & C_{p^{p},s_{t+1}^{s}}^{p} =   \sum_{n \in \mathcal{N}_{P}} \gamma_n \cdot C_{n}^{p} \label{eq:SOS2_decision} \\
        & p^{p} = \sum_{n \in \mathcal{N}_{P}} \gamma_n \cdot P_n^{p} \label{eq:SOS2_binding} \\
        & \sum_{n \in \mathcal{N}_{P}} \gamma_n = 1 \label{eq:SOS2_var} 
    \end{flalign}
\end{subequations}

\subsubsection{Energy balance and electric vehicle}
The electric energy balance of the house is given in (\ref{eq:Energy_balance}). The EV section is formulated in (\ref{eq:EV_soc_balance}) to (\ref{eq:EV_soc_lim}). The EV is modelled as a battery that has an availability pattern based on the stochastic variable $\delta^{EV}_t$. The EV discharges electricity as a load only when it has departed and can only be charged when present. The EV has a constraint in (\ref{eq:EV_soc_lim}) that states the EV must be within a certain range in its state-of-charge, and this boundary is dependent on if the EV is present, travelling, or about to travel. If it is about to depart, the SoC has a stricter boundary, to meet the safety margins of the user.
\begin{align}
    & y_t^{imp} - y_t^{exp} + y_t^{PV}  + y_t^{B,dch} = \nonumber \\
    & D_t^{El} + y_t^{EV,ch} + q_t^{sh} + y_t^{B,ch} \quad \forall t \label{eq:Energy_balance}
\end{align}
\begin{subequations}
    \begin{align}
    & SoC_t^{EV}-SoC_{t-1}^{EV} = y_t^{EV,ch}\eta^{EV}_{ch}\delta_t^{EV} \nonumber \\
    & -D^{EV} (1-\delta_t^{EV}) \quad  \quad \forall t \backslash t \neq 1 \label{eq:EV_soc_balance} \\
    %& SoC_t^{EV} - SoC^{EV,Start/Stop} = \nonumber \\
    %& y_t^{EV,ch} \eta^{EV,ch}-D^{EV}\cdot (1-\delta_t^{EV}) \quad t = 1 \label{eq:soc_balance_EV} \\
    %& SoC^{EV,Start/Stop} = SoC_t^{EV} \quad t = 24 \label{eq:soc_balance,init} \\
    & 0 \leq y_t^{EV,ch} \leq \dot{E}^{Max} \quad \forall t \label{eq:EV_charge_lim} \\
    & SoC_t^{EV,min} \leq SoC_t^{EV} \leq SoC_t^{EV,max} \quad \forall t \label{eq:EV_soc_lim} 
    \end{align}
\end{subequations}\\
In this paper, the EV energy discharge ($D^{EV}$) is considered as an input parameter per hour of unavailability. Thus, the total energy discharge varies based on the travelling duration. Introducing uncertainty to this component  would lead to more uncertain input which increases the dimension of the problem.\\

\subsubsection{Battery}
The building has a stationary battery installed that can be charged or discharged whenever needed in (\ref{eq:BAT_balance}) to (\ref{eq:Bat_soc_lim}). The battery has a specific capacity and must keep the SoC within a range to ensure the battery is not in risk of damage.
\begin{subequations}
    \begin{align}
    & SoC_t^{B} - SoC_{t-1}^{B} = y_t^{B,ch}\eta^{B}_{ch} - \frac{y_t^{B,dch}}{\eta^{b}_{dch}} \; \quad \forall t \backslash t \neq 1 \label{eq:BAT_balance} \\
    %& SoC_t^{Bat} - SoC^{bat,{Start/Stop}} = y_t^{bat,ch} - \frac{y_t^{Bat,dch}}{\eta^{bat,dch}} \quad t = 1 \label{eq:BAT_balance_init} \\
    %& SoC_t^{Bat} = SoC^{Bat,{Start/Stop}}  \quad t = 24 \label{eq:BAT_balance_end} \\
    & 0 \leq y_t^{bat,ch} \leq \dot{E}^{Bat,ch} \quad \forall t \in \mathcal{T} \label{eq:Bat_charge_lim} \\
    & 0 \leq y_t^{bat,dch} \leq \dot{E}^{bat,dch} \quad \forall t \in \mathcal{T} \label{eq:Bat_discharge_lim} \\
    %& 0 \leq y_t^{B,ch}, y_t^{B,dch} \leq \dot{E}^{B,ch},\dot{E}^{B,dch} \quad \forall t\label{eq:Bat_discharge_lim} \\
    & SoC^{B,min} \leq Soc_t^{B} \leq SoC^{B,max} \quad \forall t \label{eq:Bat_soc_lim}
    \end{align}
\end{subequations}
\subsubsection{Space heating}
All considerations regarding heating of the building is formulated in (\ref{eq:Heater}) to (\ref{eq:Envelope_balance}). The building has an electric radiator for space heating that can be turned on/off continuously. The heat dynamics in the building is seen as a grey-box model, in which the physical behavior is formulated using linear state-space models \cite{Sonderegger1978DynamicParameters, Bacher2011IdentifyingBuildings}. Using a state-space model, the dynamics between the interior temperature and the outdoor temperature can be captured alongside disturbances as heat input, which will include the impact of time-dependent temperature deviations. The heat system is represented as a 2R2C, where the building envelope is included as a middle-area between the interior and outdoor \cite{Bacher2011IdentifyingBuildings}. The only disturbances in the system considered is the heater. The interior temperature has a time-dependent boundary that must be held.
\begin{subequations}
    \begin{align}%& T_t^{in} = M^{T_{Set}} \quad t = 24 \label{eq:Heat_balance_end} \\
    & 0 \leq q_t^{sh} \leq Q^{sh} \quad \forall t \label{eq:Heater} \\
    & T_t^{in,min} \leq T_t^{in} \leq T_t^{in,max} \quad \forall t \label{eq:Heat_lim_empty} \\
    & T^{in}_{t} - T^{in}_{t-1} = \frac {1}{R_{ie}C_i} [T^{e}_{t-1} - T^{in}_{t-1}] + \frac{1}{C_i}q_{t}^{sh} \quad \forall t \backslash t \neq 1 \label{eq:interior_balance} \\
    & T^{e}_{t} - T^{e}_{t-1} = \frac {1}{R_{ie}C_e}[T^{in}_{t-1} - T^{e}_{t-1}] \nonumber \\
    & + \frac{1}{R_{ea}C_i} (T_{t-1}^{out}-T_{t-1}^e) \quad \forall t \backslash t \neq 1 \label{eq:Envelope_balance} 
    \end{align}
    \label{eq:heat}
\end{subequations}
\subsubsection{Initial conditions} \label{initial_condition}

All variables that have some energy storage property are given an initial value at the beginning/end of the scheduling day. These variables are $T_t^{in}$, $T_t^e$, $SoC_t^{EV}$ and $SoC_t^{B}$. All energy equilibrium equations have an initial equation for the time step $t=1$, being (\ref{eq:EV_soc_balance}), (\ref{eq:BAT_balance}), (\ref{eq:interior_balance}) and (\ref{eq:Envelope_balance}), where the previous time step is replaced with an initial value condition. This same is true for the last time step $t = 24$, setting the values to the initial condition. This is necessary to make it possible to represent the problem using SDP when decoupling each day, to make the end value of one day the same as the start value of the next day. This is a simplification as none of these variables are included as state variables, which is done to simplify the model.
\subsection{Solution Strategy} \label{sol_strat}
\begin{algorithm}
\For{$g = G, G-1,..,1$}
 {
    \For {$n \in N_{P}$}
    {
    $P_0^{p} \leftarrow P_n^{p}$\\
        \For{$s_g^s \in N_{S}$}
        {
            $\{ P_t^{spot}, D_t^{El}, \delta_t^{EV}, y_t^{PV},T_t^{out}, \delta_t^{Occ}\} \leftarrow StochVar(g,{s_g^s}) $
            
            $C^{p}_{i} \leftarrow \Phi(i,s_{g+1}^s,g+1)$ for $i = 1..N_{P}$
            
            $C_{s_g^s,n} \leftarrow Optimize (\ref{eq:Obj_function}) - (\ref{eq:heat})$ 
            
        }
        
        \For{$s_g^s \in N_{S}$}
        {
        $\Phi(n,s_g^s,g) = \sum_{s^s_{g+1} = 1}^{N_{S}} C_{s_g^s,n}\cdot \rho(g,s^s_{g+1}|s_g^s)$ 
        }
    }
 }
\caption{The SDP algorithm}
\label{Alg_1}
\end{algorithm}
The solution strategy for the SDP is shown in Algorithm \ref{Alg_1}. The SDP algorithm sets up the optimization problem for the building and goes backwards in the procedure to find the expected future cost curves (EFCCs) for the start of the period. We compute this for the number of discrete points $n \in N_P$ specified in line 2-3 and the number of scenarios $s_g^s \in N_S$ given in line 4. For each scenario, we realize the stochastic variables with values from StochVar in line 5, and in line 6 the EFC for the future scheduling day $g+1$ is realized in $C_i^p$. StochVar consists of the realized stochastic variables based on the scheduling day $g$ and scenario $s_g^s$. The optimization problem is then solved in line 7 to find the daily decision problem result, which is denoted as $C_{s_g^s,n}$. Furthermore, the EFC $\Phi(n,s_g^e,g)$ for each node $s_g^s$ is calculated based on the results in $C_{s_g^s,n}$ and the probability $\rho(g,s^s_{g+1}|s_g^s)$ in line 10, whereas the latter is denoted as the Markov-specific probability for the future scenario $s^s_{g+1}$ based on the current scenario $s_g^s$ and day $g$.

The combined EFCs for all $n \in \mathcal{N}_{P}$ makes up for the resulting EFCC for the given node $s_g^s$. Thus, there exists an EFCC for every scenario, which can be used in the next stage $g-1$ to represent the future uncertainty up to this point. For the initial case of $g = G$, $\Phi(n,s_g^e,G+1)$ represents the ending cost for the MP tariff for all $n \in N_P$. 

%% file: chapters/Case_study.tex
\section{Case Study} \label{Case_study}

The presented model has been applied to case study of a Norwegian building, exposed to the presented MP grid tariff. The building used is considered a single-family house (SFH) placed in the south-eastern part of Norway, and smart controls are used for the different applications. The analysis is for January 2017 with hourly time resolution per day, and the stochastic variables consist of historical data or assumed behavior. A total of four stochastic nodes per day have been generated to illustrate the uncertainty in future scheduling.

\subsection{Building structure}
The building represented in this work is modelled as a single room. For the heat dynamics of the building, the resistances and heat capacity values are based on observed values from the Living Lab building built by Zero Emission Building and NTNU \cite{LivingLab}. The interior temperature must stay between $20.5-24$ $ ^\circ C$ or $19-26$ $^\circ C$, if residents are home or not, respectively. The building is heated through a 3 $kW$ radiator. A 24 $kWh$ EV is considered with a 3.7 $kW$ smart charger, and must be between 20-90 \% capacity at all time, and between 60-90 \% when departing. When the EV is not present, a constant load of 1 $kWh$ is assumed per hour to simulate discharge due to driving. The stationary battery is based of a 5 $kWh$, 2.5 $kW$ battery from SonnenBatterie \cite{SonnenBatterie}, and must be between 10-100 \% SoC. The PV system consist of 4.65 kW installed capacity of PV with a 95 \% efficient inverter. 

The energy storage variables, as mentioned in Section \ref{Model_build}, have the following start and end condition values: $T_0^{in} = 22$ $^\circ C$, $T_0^e = 20$ $^\circ C$, $SoC_0^{EV} = 60$ $\%$, $SoC_0^{B} = 50$ $\%$.
The MP grid tariff is based on the proposal from NVE in 2017\cite{Norgesvassdrags-ogenergidirektorat2017ForslagNettvirksomhet}, adjusted for a 1 month period, with a volumetric cost of $0.00625$ $\frac{EUR}{kWh}$ and a power tariff at $7.2075$ $\frac{EUR}{\frac{kWh}{h}}$ for the month, including  25 \% VAT.

\subsection{Stochastic variables}
The stochastic variables have different scenarios between each stochastic node and day, with a total of  $N_S = 4$ scenarios per day each. Each scenario probability is only dependent on the current scenario and the future scenario, upholding the Markovian setup. The scenarios do not have the same probability; for a scenario transition on the same node, the probability is higher than the others. The electricity prices are obtained from NordPool, and only January 2017 for NO1 has been used. The PV irradiation and outdoor temperature series are from historical data for January 2014-2017 obtained from Nibio \cite{NedlastingLMT}, while the electric-specific non-shiftable load are from different buildings in January 2017, but with similar consumption obtained from Ringerikskraft \cite{Ringerikskraft}, a power company in southern Norway. The occupancy for EV, which also affects temperature boundary, is based on assumed behaviour regarding leaving for work and weekend travelling, with four different patterns.

\subsection{Model cases}
To keep the level of accuracy for the EFCC piecewise-linear function, the peak power value was discretized to be $\mathcal{N}_{P} = 41$, ranging from 0-10 $\frac{kWh}{h}$. This gives a total of $5084$ nodes to analyze for the SDP algorithm.
The main case $Peak_{SDP}$ will analyze how the proposed SDP algorithm contributes with finding the cost-optimal initial peak power for the building with the MP grid tariff, and what the expected total cost for the residents will be over the month. To achieve this, the EFCC will be used in a simulation phase where the days are run sequentially day by day, with the peak power achieved included in the transition and the EFCC representing the future, to find the total cost for that simulated month. The transition of each day is decided by the Markovian probability. 1000 simulations will be performed to get an accurate description of the uncertainty in the model.

In addition to the main case, there will be 4 other cases analyzed and compared against. Case $Peak_{no}$ is a case where the resident does not consider the peak power in the scheduling decision and case $Peak_{min}$ is where the resident minimizes the peak from day one without any knowledge of the future. The last two cases considers the entire month scheduled as one holistic problem, with perfect information on the whole period. Case $Hol_{init}$ is with the initial condition for the energy variables being kept true for every 24 hours, while case $Hol$ is without the initial condition except for the last scheduling day.

%% file: chapters/Results_and_discussions.tex
\section{Results \& Discussion} \label{discussion}

\subsection{Expected future cost curves}
\begin{figure}[h!]
\centering
\includegraphics[width=0.50\textwidth]{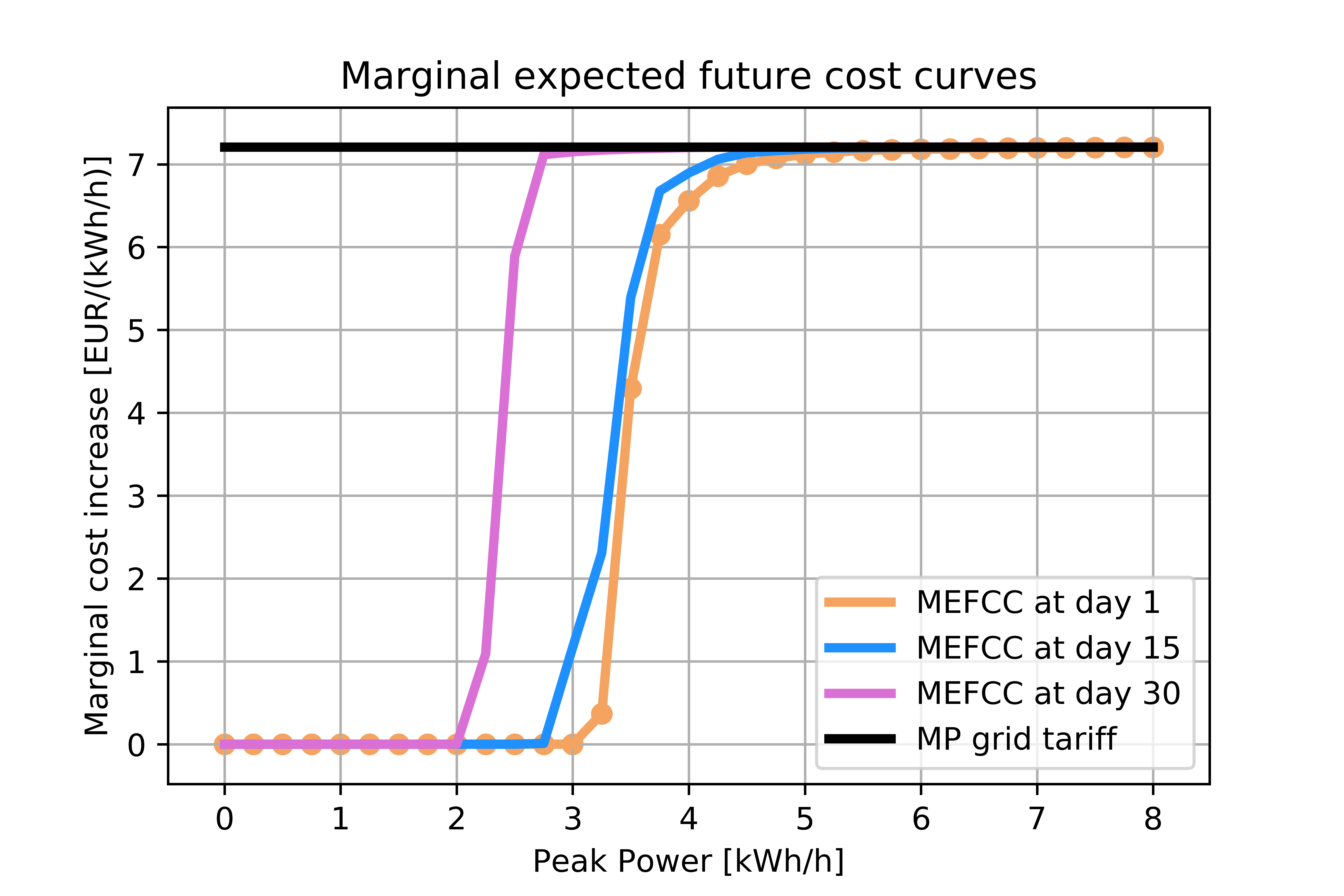}
\caption{Plot of the non-linear marginal expected future cost curve for three different days and a given scenario node.}
\label{fig:EFCC}
\centering
\end{figure}

From the SDP algorithm in Algorithm \ref{Alg_1}, EFCCs are created for every day of the month, and for every scenario within the days, based on the highest peak power. For each additional future day the EFCC represents, the total future cost will increase. To give a comparable overview of the EFFCs over the period, the values plotted in Fig. \ref{fig:EFCC} are the marginal EFCC of increasing the peak power by 1 $\frac{kWh}{h}$.

In general, each EFCC has three main parts that can be distinguished. The start of the curves are flat, where initializing on these values will result in a re-adjustment later on. The middle section of the curves are where the non-linearity is present, as increasing the initial peak power will result in a marginal cost increase. The varying marginal value is a trade-off between paying for more flexibility from import versus scheduling with less flexibility. The end of the curve is flat as the increase in import peak gives no extra incentives in scheduling and only a cost increase equal to the MP grid tariff.

% and are non-linear cost curves based on the highest peak power.
%which has an impact on the curve and depicts the starting flat curve level for all days prior.
%Note that this does not mean that the EFCC gives the perfect level to aim for that will always be held. 
%the weighted future cost based on the peak power. As this methodology is based on scenarios and probabilities, the curves and their proposed peak power level is the cost-optimal based on the uncertainty in the system, where the impact each scenario has is reflected with their cost and probability. 
%to benefit from this every sequential day later. if the savings of having more import capability exceeds the EFC at the current peak power, it is more optimal to increase the peak capacity as this will induce more savings in the future.

The plot shows that the point in which the non-linear curve starts appearing changes as more days are included. If a day has scenarios with high demand, the optimal peak power might have to be larger than calculated in later days, which changes the starting flat curve level further on. Thus, the critical days will be reflected in the EFCC for day 1. However, as uncertainty is included in the model, these EFCC only gives the cost-optimal future cost based on probability. Thus, the cost-optimal peak power from the EFCC can give simulations where the peak power is exceeded due to critical scenarios being realized, and simulations where the actual lowest peak power level would be below the initial value. This is the consequence of adding uncertainty into the mix. 

Based on the marginal EFCC in Fig. \ref{fig:EFCC}, the optimal initial level of the peak power in day 1 is 3 $\frac{kWh}{h}$, as that is the highest peak power with the lowest cost. For each day, the EFCC curves are included to reflect the future impact on the current decision-making, and indicates what cost-optimal level of peak power the building should try to achieve. It considers the cost of load shifting to reduce the peak power to the cost of increasing peak power and the future impact, and optimizes based on the expected cost-effective solution.% Thus, the model can decide to both shift the load to more expensive hours, or increase the peak power if the load shifting is more expensive.

%For each day, if the savings within the day is higher than the cost for extra peak power capacity, then it will be more beneficial to increase the capacity as this benefit will continue in the future. However, there will be scenarios where this is not necessary and will lead to a higher total cost, and thus the scenario and the value of flexibility will affect the optimal peak capacity.

\subsection{Simulation results}

The results from the simulations regarding all cases are shown in Table \ref{Table_Results}, with both the average total cost, and the average peak power achieved. The SDP algorithm manages to keep the average total cost below what would be the case if the resident was unaware and had no prediction of the future scenarios, when comparing case $P_{SDP}$ to cases $P_{no}$ and $P_{min}$. %For each day, the EFCC curves are included to reflect the future impact on the current decision-making, and indicates what cost-optimal level of peak power the building should try to achieve. The internal flexibility within space heating, water tank, EV charging and battery usage enables load shifting to contribute with lowering the peak power. It considers the cost of load shifting to reduce the peak power to the cost of increasing peak power and the future impact, and optimizes based on the cost-effective solution. Thus, the model can decide to both shift the load to more expensive hours, or increase the peak power as the load shifting is more expensive.

Case $P_{no}$ uses the smart control and the internal flexibility to minimize cost from the variation in spot prices, ignoring the penalty paid at the end of the month in the daily decision-making. Thus, the cost increase of 36.1 \%, compared to the $P_{SDP}$ result, is to be expected as the achieved peak is $10 \frac{kWh}{h}$. For case $P_{min}$, the user has no information on the future predictions, and react by keeping the peak as low as possible by considering achieved values. The loss of value due to the effort of load shifting to keep the import level at an unnecessary low level early, is the reason this setup is 0.3 \% higher in average total cost, despite having a 0.3 \% lower average import peak, compared to $P_{SDP}$.

\begin{table}[htbp]
\centering
\caption{Average total cost and peak power for the 5 cases.}
\begin{tabularx}{\linewidth}{X*{9}{c}}
\hline
\textbf{Cases} & $P_{SPD}$ & $P_{no}$ & $P_{min}$ & $P_{Hol}$ & $P_{Hol,init}$ \\ \hline
Total cost [EUR] & 131.2 & 178.6 & 131.6 & 126.5 & 130.4 \\ \hline
Peak [$\frac{kWh}{h}$] & 3.24 & 10.0 & 3.23 & 2.69 & 3.22 \\ \hline
\end{tabularx}
\label{Table_Results}
\end{table}

%This case will then utilize the internal flexibility to shift loads to keep the import profile flat, but at a level that will be exceeded on later days. 

A drawback of the SDP algorithm is that for each day, the start and end condition for variables considering energy storage must be equal, as mentioned in \ref{initial_condition}. This is not the optimal case if the whole period is already known, due to load shifting between days. For the holistic case with the initial condition limitations $P_{Hol,init}$ kept, the difference of the average total cost is lowered by 0.6 \% and 0.6 \% lower average peak power compared to $P_{SDP}$. However, when neglecting the limitation in case $P_{Hol}$, the average total cost is 3.6 \% lower compared to $P_{SDP}$, and the average peak power is 17.0 \% lower. The results from both cases with perfect information illustrates the weakness of simplifying energy storage values when transitioning from one day to the next, as this increases the total cost. The current SDP algorithm cannot load shift between the days, limiting the complete benefit of flexibility which is shown in the holistic cases. This shows the value of deviating from the initial conditions. % the resident has limited benefits if the initial conditions must be held. %If the resident had perfect information, but restrained to keep the initial condition through every day, then the benefit from the perfect information is small compared to if he could neglect this constraint.
%This impact is shown in the holistic solutions.

We demonstrated  the benefit of including not only the peak power from the grid, but other information of the building that is carried over between scheduling days. This is a concept that the SDP algorithm can include, however, this can create scalability issues if the number of state variables increases too much. Including the critical energy variables and simplifying the rest can improve the performance. Another point of interest is to find logical and season-optimal initial conditions for the energy variables. The ones for this case study might be ill-suited for the transition. The findings from the holistic solutions support that the SDP algorithm can improve in performance if the initial conditions are analyzed further, both regarding initial values and including them as state variables in the algorithm.

%% file: chapters/Conclusion.tex
\section{Conclusion} \label{Conclusion}

With new capacity based grid tariffs being proposed in Norway by NVE, the value of load shifting for the end-users will come in focus. If the total cost depends on a longer horizon, the long-term impact should be represented today, especially with uncertainty included. We present a model that aims at representing the future cost for the system based on current decision-making for a building when considering the MP grid tariff. The model manages to represent the EFC for the building as a cost curve depicting the optimal peak power to aim for based on expected future behaviour. The Markovian setup defines the probable future states of the problem, which depend only upon the present state, and are not conditioned on the sequence of states and actions that preceded them. If these dependencies are considered, it will cause an explosion in the size of the state representation, and correspondingly, the algorithm becomes computationally infeasible.

%finding the optimal highest peak power that is cost-effective for a building based on uncertainty in load demand, user behavior, external impact such as irradiation and temperature, by utilizing smart control of applications for load shifting. The model manages to represent the expected future cost for the building based on highest peak power as a non-linear cost curve, depicting the optimal peak power to aim for based on the previously achieved level, current decision-making and expected future behaviour. 

This model was applied to a realistic Norwegian household for January 2017 to find the optimal peak power and the expected total electricity cost. This was compared to two situations where the future cost from the grid tariff is unknown or not considered for the building, and two cases with perfect information regarding stochastic inputs. For the cases where the peak power is either minimized daily or not considered into the daily scheduling, they have a 0.3 \% and 36.1 \% cost increase compared to the proposed method, respectively. For holistic setups where the initial condition was kept at each day or ignored, the results showed a cost save of 0.6 \% and 3.6 \% compared to the proposed methodology, respectively. This showcases the importance of these initial conditions.

%Question: Reviewer 3 noted that the Markov decision should be emphasized in the conclusion about the property assumption. My comment is that i am unsure about what he means:
%       I am a bit unsure about how I should answer this. Is he questioning my values, or is he asking me to specify what value this could add to the scheduling? Is it about why we use Markov, or is it about what other methods could be applied?

%The proposed model finds the expected future cost based on peak power for the start of the period by analyzing the future stages and the uncertainty and finds the cost-optimal initial peak power for the period.
% However, the results also indicated that the model is limited by the use of initial conditions for several variables which must be held at the beginning and end of each day.
%When compared to a holistic solution with perfect information for the whole month, but with the initial conditions kept, the holistic solution had a cost save of 0.6 \& compared to the proposed methodology. A holistic solution without the initial condition limitation had a cost reduction of 4.2 \%, showcasing the importance of these initial conditions. 
%The results showed that this representation reduced the electricity cost by 0.2 \% compared a minimization of the peak based on current and prevision decision-making. If the peak power was ignored when scheduling, then the expected cost would be 26.6 \% higher than the proposed method. 

\section{Future Work}

The findings showcase the potential the representation of the long-term uncertainty has on the overall result. %The curves could, if accurately describing the future uncertain possibilities, be connected to a detailed short-term operational model to represent the far end of the optimization horizon.
However, the results also indicated that further investigation into the initial conditions that consider energy storage in the building must be done to find the impact they have on the system. Whether they need to be given an accurate initial condition or included into the future cost curve will be a key question to answer.

\section{Acknowledgment} \label{Acknowledgment}

This work was funded and supported by the Research Council of Norway (Grant Number: 257626/257660/E20) and several partners through FME ZEN and FME CINELDI. The authors want to thank Karen Byskov Lindberg from SINTEF Community and her team with assistance in formulating the state-space model for space heating. Also, we thank Ringerikskraft for giving electricity load profiles from a selection of their households.